\documentclass[a4paper,twoside]{article}
\usepackage{xypic}
\usepackage{babel}
\usepackage{amssymb}

\def\R{{\mathbb R}} 
\def\C{{\mathbb C}}

\begin{document}
\baselineskip=0.45cm

\thispagestyle{empty}

\pagestyle{myheadings}

\markboth{{\small { \'Elimination des singularit\'es pour les 
fonctions CR}}}{{\rm J. {\sc Merker} and E. {\sc Porten}}}

$\:$
\bigskip

\noindent
{\Large {\bf Enveloppe d'holomorphie locale des
vari\'et\'es CR 

\bigskip

\noindent
et \'elimination
des singularit\'es 

\bigskip

\noindent
pour les fonctions CR int\'egrables}}

\bigskip

\noindent
{\large Projet de note, par {\bf Jo\"el MERKER${\:\!}^{*}$} 
et {\bf Egmont PORTEN${\:\!}^{**}$}}.

\bigskip

\noindent
$(*)$ LATP, CMI, 39 rue Joliot Curie, F-13453 Marseille Cedex 13\\
Courriel : merker@gyptis.univ-mrs.fr

\smallskip

\noindent
$(**)$ Max-Planck-Gesellschaft, Humboldt-Universit\"at zu Berlin\\
J\"agerstrasse, 10-11, D-10117 Berlin, Germany\\
Courriel : egmont@mathematik.hu-berlin.de

\medskip

{\bf R\'esum\'e.}
\begin{minipage}[t]{11.75cm}{\small 
Soient $M$ une vari\'et\'e CR localement
plongeable et $\Phi\subset M$ un ferm\'e. On donne des conditions
suffisantes pour que les fonctions $L_{loc}^1$ qui sont 
CR sur $M\backslash \Phi$ le soient aussi sur $M$ tout entier.}
\end{minipage}

\bigskip

\noindent
\hspace{2.5cm}{\bf Local envelope of holomorphy
of CR manifolds

\noindent
\hspace{2.5cm} and removable singularities for integrable CR functions.}

\bigskip

{\bf Abstract.} 
\begin{minipage}[t]{11.75cm}{\small 
Let $M$ be a locally embeddable CR manifold
and $\Phi\subset M$ be a closed set. We 
give sufficient conditions in order 
that $L^{1}_{loc}$ functions on $M$ which are CR
on $M\backslash \Phi$ are CR on $M$. }
\end{minipage}

\bigskip
\bigskip
\bigskip

\noindent
{\bf Abridged English Version.} Let $M$ be a locally embeddable
CR manifold, ${\rm dim}_{CR} M=m$, 
${\rm codim} \ M=n$, ${\rm dim} \ M=d=2m+n$ and
$\Phi\subset M$ a closed set.
We give various conditions in order that $\Phi$ is $L^1$-removable, {\it i.e.}
                      \begin{equation}
                      L^1_{loc}(M) \cap L_{loc,CR}^1(M\backslash\Phi
)=L_{loc,CR}^1(M)
                      \end{equation} Let $H^{\kappa}$ denote
$\kappa$-dimensional Hausdorff measure.

\smallskip \noindent {\sc Theorem 1.} -- {\it If $M$ is ${\cal C}^3$,
a function $f\in L_{loc}^1(M)$ is CR if and only if $f|_{{\cal O}}$
belongs to $L_{loc,CR}^1({\cal O})$ for almost every CR orbit ${\cal
O}$.}

\smallskip \noindent {\sc Corollary 1.} -- {\it If $\Phi=\cup_{a\in A}
{\cal O}_a$ is of zero $d$-dimensional measure, then $(1)$ is
satisfied.}

\smallskip Theorem 1 reduces the problem to the case where $M$ is a
single CR orbit, {\it i.e.} $M$ is {\it globally minimal} \cite{ME1}.

\smallskip \noindent {\sc Theorem 2.} -- {\it Let $M$ be ${\cal
C}^{2,\alpha}$, $0\leq \alpha <1$, ${\rm dim}_{CR} M =m\geq 1$. Every
closed subset $E$ of $M$ such that $M$ and $M\backslash E$ are globally
minimal and such that $H_{loc}^{d-3}(E)<\infty$ is $L^1$-removable.}

\smallskip The notion of wedge (${\cal W}$-)-removability is defined
here in higher codimension.

\smallskip \noindent {\sc Theorem 3.} -- {\it Let $M$ be ${\cal
C}^{\omega}$, ${\rm dim}_{CR}M=m\geq 1$. Every closed set $E\subset M$
such that $M$ and $M\backslash E$ are globally minimal and such that
$H^{d-2}(E)=0$ is ${\cal W}$- and $L^1$-removable.}

\smallskip \noindent {\sc Theorem 4.} -- {\it Let $M$ be ${\cal
C}^{2,\alpha}$, $m\geq 1$, and let $N$ be a ${\cal C}^2$ connected
submanifold of $M$ such that $M$ and $M\backslash N$ are
globally minimal.

{\rm (i)} \ If ${\rm codim}_M N \geq 3$, then $N$ is ${\cal
W}$- and $L^1$-removable{\rm ;}

{\rm (ii)} \ Every closed set $\Phi \subset N$ is ${\cal
W}$- and $L^1$-removable if $\Phi \neq N$, ${\rm codim}_M N =2$ and $m\geq
1${\rm ;}

{\rm (iii)} \ $N$ is ${\cal W}$- and $L^1$-removable if $N$ is generic at
one point, ${\rm codim}_M N =2$ and $m\geq 2$.}

\smallskip A set $S\subset M$ is called a ${\cal C}^{\lambda}$ peak
set, $0<\lambda <1$, if there exists a {\it nonconstant} function $\varpi
\in {\cal C}^{\lambda}_{CR}(M)$ such that $S=\{\varpi=1\}$ and
$|\varpi| \leq 1$.

\smallskip \noindent {\sc Theorem 5.} -- {\it Let $M$ be ${\cal
C}^{2,\alpha}$ globally minimal. Then every ${\cal C}^{\lambda}$ peak
set $S$ satisfies $H^{d}(S)=0$ and is $L^1$-removable.}

\smallskip \noindent {\sc Corollary 2.} -- {\it Let $M$ be ${\cal
C}^3$. Then a ${\cal C}^{\lambda}$ peak set $S$ is $L^1$-removable if
$H^d(\cup_{{\cal O}\subset S} {\cal O})=0$.}

\bigskip

\def\theequation{\arabic{equation}}
\setcounter{equation}{0}

\noindent
{\large {\bf 1. \'Enonc\'es}}

\medskip

Soient $M$ une vari\'et\'e CR localement
plongeable, de dimension CR ${\rm dim}_{CR} M =m$, 
de codimension ${\rm codim} \ M= n$, 
de dimension ${\rm dim} \ M= 
d=2m+n$
et soit $\Phi\subset M$ un 
ferm\'e de $M$. 
Dans ce travail,
on cherche des conditions, portant sur $M$ et $\Phi$, pour que l'on ait
                      \begin{equation}
                      L_{loc}^1(M)\cap L_{loc,CR}^1(M
                      \backslash \Phi)=L_{loc}^1(M).
                      \end{equation}
Si (1) est v\'erifi\'ee, on dira que $\Phi$ est $L^1$-\'eliminable.
D'apr\`es Tr\'epreau \cite{TR}, 
toute vari\'et\'e CR $M$ est r\'eunion disjointe
(en g\'en\'eral transfinie)
de sous-vari\'et\'es CR immerg\'ees connexes 
$M=\cup_{i \in {\cal I}} {\cal O}_i$,
${\cal O}_i \subset M$, appel\'ees {\it orbites CR de $M$}, qui sont 
{\it caract\'eristiques}, {\it i.e.} ${\rm dim}_{CR} {\cal O}_i={\rm dim}_{CR} M$ et
minimales pour l'inclusion et cette propri\'et\'e. Un premier 
r\'esultat concerne ces vari\'et\'es.
Il a \'et\'e d\'emontr\'e par B. J\"oricke dans la classe
${\cal C}^0(M)$ \cite{JO1} et dans 
$L_{loc}^1(M)$ si $M$ est 
une hypersurface de classe ${\cal C}^2$ \cite{JO2}. Enfin, le second
auteur l'a \'etendu \`a $L_{loc}^1(M)$ en codimension 
quelconque dans sa th\`ese \cite{PO}.

\smallskip
\noindent
{\sc Th\'eor\`eme 1.} -- {\it Si $M$ est de classe
${\cal C}^3$, une fonction $f\in L_{loc}^1(M)$ est CR si et seulement
si $f|_{{\cal O}_i}$ appartient \`a $L_{loc,CR}^1({\cal O}_i)$ pour presque 
toute orbite CR ${\cal O}_i$, au sens de la mesure sur $M$.}

\smallskip
La restriction $f|_{{\cal O}_i}$ est bien d\'efinie et appartient
\`a $L_{loc}^1({\cal O}_i)$ pour presque toute orbite ${\cal O}_i$.
Nous renvoyons le lecteur \`a 
\cite{PO} ou \cite{MP1} pour une preuve compl\`ete.

\smallskip
\noindent
{\sc Corollaire 1.} -- {\it Si $\Phi=\cup_{a\in A} {\cal O}_a$ est de
mesure $d$-dimensionnelle nulle, alors (1) est v\'erifi\'ee.}

\smallskip

Le th\'eor\`eme 1 r\'eduit l'\'etude de (1) au cas
o\`u $M={\cal O}_i$ est une seule orbite. 
L'aspect central de notre travail consiste justement 
\`a replacer l'\'etude de (1) dans le contexte de la th\'eorie des
orbites CR et de l'extension des fonctions CR, b\'en\'eficiant
en cela des travaux de Tr\'epreau, Tumanov et J\"oricke. Lorsque
$M$ est une hypersurface, le probl\`eme (1) est trait\'e par
J\"oricke et Chirka-Stout \cite{CS}.

Bien entendu,
les orbites sont aussi localement plongeables. Soit
$M$ g\'en\'erique dans $\C^{m+n}$. Dans ce cas, un ouvert
connexe ${\cal W}_0$ sera appel\'e {\it wedge attach\'e}
\`a $M\backslash \Phi$ s'il existe une section continue
$\eta: M\backslash \Phi \to T_M\C^{m+n}\backslash \{0\}$
du fibr\'e normal \`a $M$ telle que
${\cal W}_0$ contient un wedge ${\cal W}_p$
d'edge $M$ en $(p,\eta(p))$ ({\it cf.} \cite{TU}, p.3), pour tout
point $p\in M\backslash \Phi$. Cette notion a un sens local
lorsque $M$ est localement plongeable. Le probl\`eme (1) fait 
intervenir la g\'eom\'etrie des wedges attach\'es.

\smallskip
\noindent
{\sc D\'efinition 1.} $\Phi$ est dit ${\cal W}$-\'eliminable
si, pour tout wedge ${\cal W}_0$ attach\'e \`a 
$M\backslash \Phi$, il existe un wedge 
${\cal W}$ attach\'e \`a $M$ tel que les fonctions 
holomorphes dans ${\cal W}_0$ se prolongent
holomorphiquement \`a ${\cal W}$. 

\smallskip
On note $H^{\kappa}(E)$ la mesure de Hausdorff $\kappa$-dimensionnelle
de $E$, pour une m\'etrique fix\'ee sur $M$.
$H^d$ s'identifie \`a la mesure de Lebesgue. 
Enfin, on dira qu'une vari\'et\'e CR est {\it globalement
minimale} si elle consiste en une seule orbite CR. Le ferm\'e $\Phi$
sera not\'e $E, N$ ou $S$, suivant le contexte.

\smallskip
\noindent
{\sc Th\'eor\`eme 2.} -- {\it Soit $M$ de classe
${\cal C}^{2,\alpha}$, $0<\alpha <1$, ${\rm dim}_{CR} M=m\geq 1$.
Tout ferm\'e $E$ de $M$ tel que $M$ et $M\backslash E$
sont globalement minimales et tel que $H_{loc}^{d-3}(E)<\infty$
est $L^1$-\'eliminable.}

\smallskip
Par exemple, sous ces hypoth\`eses, toute sous-vari\'et\'e
$N$ de codimension au moins trois, est $L^1$-\'eliminable.
Le r\'esultat suivant (\cite{MP2})
s'applique \`a l'extension des fonctions CR m\'eromorphes.

\smallskip
\noindent
{\sc Th\'eor\`eme 3.} -- {\it Soit $M$ de classe ${\cal C}^{\omega}$,
${\rm dim}_{CR}M=  m\geq 1$. Tout ferm\'e $E\subset M$
tel que $M$ et $M\backslash E$ sont globalement
minimales et tel que $H^{d-2}(E)=0$ est
${\cal W}$- et $L^1$-\'eliminable.}

\smallskip
Le cas de singularit\'es plus massives est trait\'e dans
\cite{MP1}:

\smallskip
\noindent
{\sc Th\'eor\`eme 4.} -- {\it Soient $M$ ${\cal C}^{2,\alpha}$, $m\geq 1$. Soit
$N$ une sous-vari\'et\'e connexe de $M$, de classe ${\cal C}^2$
telle que $M$ et $M\backslash N$ sont globalement minimales. Alors

{\rm (i)} \ $N$ est ${\cal W}$- et $L^1$-\'eliminable si 
${\rm codim}_M N \geq 3$.{\rm ;}

{\rm (ii)} \ Tout ferm\'e $\Phi$ de $N$ est  
${\cal W}$- et $L^1$-\'eliminable, si $\Phi\neq N$, ${\rm codim}_M N =2$
et $m\geq 1${\rm ;}

{\rm (iii)} \ $N$ est ${\cal W}$- et $L^1$-\'eliminable 
si $N$ est g\'en\'erique
en au moins un point, ${\rm codim}_M N =2$ et $m\geq 2$.}

\smallskip
{\it Remarque.} L'\'elimination $L^1$ de compacts $K\subset\subset N$ 
de vari\'et\'es $N\subset M$
g\'en\'eriques de codimension un 
appara\^{\i}t dans les travaux de B. J\"oricke
\cite{JO2} pour $n=1$ et dans \cite{PO} pour $n\geq 2$.

\smallskip
Une version plus faible du th\'eor\`eme 4 est contenue dans 
\cite{ME2} dans le cas de l'\'elimination ${\cal W}$.
L'hypoth\`ese d'orbite sur $M$ et $M\backslash N$
est essentiellement n\'ecessaire : si elle n'est pas satisfaite, il 
existe des exemples simples de $M$, $N$ et de
distributions CR de support une sous-vari\'et\'e caract\'eristique 
ferm\'ee $S$ de $M\backslash  N$ qui ne se prolongent
pas holomorphiquement \`a un wedge au-dessus de $N$. Enfin,
les hypoth\`eses g\'eom\'etriques sur $N$ et $\Phi$ sont 
calqu\'ees sur celles qui rendent les th\'eor\`emes
2,3 et 4 connus lorsque $M$ est un ouvert de $\C^m$, 
{\it i.e.}  $n=0$.

\smallskip
Maintenant, un sous-ensemble
$S$ de $M$ est dit {\it ensemble pic h\"old\'erien} s'il a la forme
$\{\varpi=1\}$, avec $\varpi\in {\cal C}^{\gamma}_{CR}(M)$ 
{\it non constante} pour un $\gamma$, 
$0<\gamma<1$ et $|\varpi|\leq 1$. Gr\^ace au th\'eor\`eme 1
et aux techniques de d\'eformation de disques,
on g\'en\'eralise les
r\'esultats de \cite{KR}.

\smallskip
\noindent
{\sc Th\'eor\`eme 5.} -- {\it Soient $M$ ${\cal C}^{2,\alpha}$ globalement
minimale. Tout ensemble pic h\"old\'erien $S$
de $M$ v\'erifie $H^d(S)=0$
et est $L^1$-\'eliminable.}

\smallskip
\noindent
{\sc Corollaire 2.}-- {\it Soit $M$ ${\cal C}^3$. Un ensemble
pic h\"old\'erien $S$ est $L^1$-\'eliminable si
$H^d(\bigcup_{{\cal O}_{\i} \subset S} {\cal O}_{\i})=0.$}

\smallskip
Enfin, puisque $L^{\rm p}_{loc}$ se plonge dans 
$L_{loc}^1$ pour ${\rm p} \geq 1$, tous ces 
r\'esultats sont valables dans $L^{\rm p}_{loc}$.

\bigskip

\noindent
{\large {\bf 2. Preuves.}}

\medskip

Nous allons donner ici un r\'esum\'e 
des preuves des th\'eor\`emes 2,3,4 et 5. Les
preuves rigoureuses sont contenues dans \cite{MP1}, \cite{MP2}, \cite{PO}.
Le th\'eor\`eme 1 pour $f\in {\cal C}^0_{CR}(M)$ est d\'emontr\'e
dans \cite{JO1} et dans \cite{PO} pour $f\in L_{loc,CR}^{\rm p}$.

C'est B. J\"oricke qui a eu l'id\'ee d'utiliser l'in\'egalit\'e
de Carleson sur des familles r\'eguli\`eres de disques
analytiques attach\'es \`a $M$ pour d\'eduire l'\'elimination
$L^1$ de l'\'elimination ${\cal W}$, dans le 
cas hypersurface. Dans ce travail et dans \cite{MP1}, \cite{MP2},
nous raffinons les r\'esultats de \cite{ME2} pour 
la ${\cal W}$-\'elimination et les \'etendons \`a $L^1$ comme
dans \cite{JO2}. Mentionnons enfin que la technique
dite de <<balayage par des wedges>> utilis\'ee dans \cite{JO2},
\cite{CS}, \cite{PO} ne s'appliquerait
qu'en dimension CR $m\geq 2$; c'est pourquoi nous
utilisons ici les d\'eformations de disques analytiques et le 
principe de continuit\'e.

\smallskip
On traite le cas $L^1$,
le cas ${\cal W}$ sera d\'emontr\'e en cours.  La technique
consiste en un grand nombre de d\'eformations de $M$ dans des ouverts
obtenus en attachant des disques analytiques \`a $M$ et \`a ses
d\'eformations. En particulier, nous d\'ecrivons une partie de
l'enveloppe d'holomorphie d'ouverts de type wedge attach\'es \`a
$M\backslash \Phi$, assez \'etendue pour proc\'eder ensuite \`a
l'\'elimination $L^1$.

\noindent {\it \'Etape 1.}  Gr\^ace au th\'eor\`eme d'extension de
Tr\'epreau-Tumanov g\'en\'eralis\'e (\cite{ME1}, \cite{ME2}) et au
th\'eor\`eme de l'<<edge of the wedge>>, on peut prolonger
holomorphiquement $f\in L_{loc,CR}^1(M\backslash \Phi)$ \`a un wedge
${\cal W}_0$ attach\'e \`a $M\backslash \Phi$. Pour le contr\^ole en
norme $L^1$ de l'extension, on utilise de <<bonnes>> familles de
disques analytiques attach\'es \`a $M$. Soit $\Delta$ le disque
unit\'e dans $\C$, $b\Delta$ son bord.

\smallskip \noindent {\sc D\'efinition 2.} -- On appelle {\it famille
r\'eguli\`ere en $p$ de disques analytiques attach\'es \`a $M$} une
application ${\cal C}^{2,\beta}$, $\beta<\alpha$, $A: {\cal S} \times
{\cal V}\times \overline{\Delta} \to \C^{m+n}$, $(s,v,\zeta)\mapsto
A_{s,v}(\zeta)$, holomorphe en $\zeta$, o\`u $(0\in)
{\cal S} \subset \R^{2m+n-1}$, $(0\in) {\cal V} \subset\R^{n-1}$ sont des
ouverts, telle que $A_{0,v}(1)=p$ et que

1) L'application ${\cal S} \times b\Delta \to M$, $(s,\zeta)\mapsto
A_{s,v}(\zeta)$ est un plongement, $\forall \ v\in {\cal V}$;

2) Le vecteur $\eta:=-\partial A_{0,0}
/\partial\zeta (1) \not\in T_pM$
et ${\rm rang}(v\mapsto {\rm
pr}_{T_p\C^{m+n}/(T_pM\oplus\R\eta)}(-\partial A_{0,v} /\partial\zeta
(1)))=n-1$.

\smallskip Toute famille r\'eguli\`ere d\'efinit un wedge ${\cal
W}_{A,p}$ en $(p,\eta)$ par ${\cal W}_{A,p}:= \{A_{s,v}(\zeta)\in
\C^{m+n}: \ (s,v,\zeta)\in {\cal S}_1\times {\cal V}_1\times
\stackrel{\circ}{\Delta}_1\}$, o\`u ${\cal S}_1 \subset {\cal S}$, ${\cal V}_1 \subset
{\cal V}$, $\Delta_1=\overline{\Delta}\cap \Delta(1,\rho_1)$,
$\rho_1>0$.  Soit ${\cal W}={\cal W}(U,C)=\{z+\eta: \ z\in U, \eta\in
C\}$ un wedge de base $U\subset M$ et de c\^one $C\subset \C^{m+n}$,
{\it e.g.} ${\cal W}\approx {\cal W}_{A,p}$.

\smallskip \noindent {\sc D\'efinition 3.} -- Une fonction $f\in
L_{loc,CR}^1(M)$ est dite {\it prolongeable dans} ${\rm H}_{\rm
a}^1({\cal W})$ s'il existe $F\in {\cal H}({\cal W})$ telle que
$F|_{U_{\eta}}\to f$ au sens $L^1$, o\`u $U_{\eta}= U+\eta$, $\eta\in
C$, uniform\'ement lorsque $\eta \to 0$.

\smallskip

Pour d\'emontrer l'extension de $f$ dans la classe de Hardy ${\rm
H}_{\rm a}^1({\cal W}_{A,p})$, on utilise des familles r\'eguli\`eres
attach\'ees \`a des d\'eformations de $M$:

\smallskip \noindent {\sc Proposition 1} (\cite{PO}). 
-- {\it Soit $M$ globalement
minimale, ${\cal C}^{2,\alpha}$.  Pour tout $\varepsilon>0$, il existe
une d\'eformation ${\cal C}^{2,\beta}$ $(d,M^d)$ 
de $M$ \`a support compact
avec $M^d\equiv M$ pr\`es de $p$ et $||M^d-M||_{{\cal
C}^{2,\beta}}<\varepsilon$, telle que

1) Il existe une famille r\'eguli\`ere de disques $A_{s,v}$ attach\'es
\`a $M^d${\rm ;}

2) Il existe un op\'erateur lin\'eaire born\'e de prolongement
$L_{loc,CR}^1(M)\to L^1_{loc,CR}(M^d)$, $f\mapsto f^d$, tel que
$f^d\equiv f$ sur l'ensemble o\`u $M^d$ co\"{\i}ncide avec $M${\rm ;}

3) Pour $f\in L_{loc,CR}^1(M)$ fix\'ee, il existe une d\'eformation
$d$ telle que $||f^d-f||_{L^1}<\varepsilon$.}

\smallskip
Cette construction demande l'emploi en d\'etail des
techniques de d\'eformation de disques de Bishop
\'elabor\'ees par Tumanov (\cite{TU})
et l'existence de $M^d$ satisfaisant $1)$ et $2)$ 
\'equivaut \`a la globalit\'e minimale (\cite{PO}, \cite{MP1}).

\smallskip
\noindent
{\sc Proposition 2.} -- {\it $L_{loc,CR}^1(M)$ est prolongeable dans
${\rm H}_{\rm a}^1({\cal W}_{A,p})$.}

\smallskip En effet, sur chaque disque remplissant
${\cal W}_{A,p}$, on se ram\`ene \`a l'estim\'ee de
Carleson sur $\Delta$: si $r(\theta)\in {\cal C}^1([-\pi,\pi],[0,1])$
telle que $r\equiv r_1$ sur $(-\theta_1,\theta_1)$, $0<r_1<1$ et
$\hbox{supp} \ (1-r)\subset 
(-\theta_0,\theta_0)$, $0<\theta_1 <\theta_0 <\pi$,
alors il existe $C>0$ telle que
$\forall \ u \in H^1_{\rm a}(\Delta), \ 
\int_{-\theta_0}^{\theta_0} 
|u(re^{i\theta})| d\theta \leq C ||u||_{L^1(b\Delta)}$. $\square$

\smallskip
\noindent
{\it \'Etape 2.}
On note ${\cal H}({\cal U})$ l'anneau des fonctions
holomorphes dans ${\cal U}$ et 
${\cal V}(E)$ un voisinage ouvert arbitrairement
petit d'un ensemble $E$ dans $\C^{m+n}$.
La d\'eformation suivante r\'eduit la d\'emarche 
au cas o\`u $L^1_{loc,CR}(M\backslash \Phi)\cap L_{loc}^1(M)$ a \'et\'e
remplac\'e par ${\cal H}({\cal V}(M^d\backslash \Phi))\cap 
L_{loc}^1(M^d)$.

\smallskip
\noindent
{\sc Proposition 3.} -- {\it Soient $M$ ${\cal C}^{2,\alpha}$, 
g\'en\'erique dans $\C^{m+n}$,
$f\in L_{loc,CR}^1(M)$ et $U$ un ouvert de $M\backslash \Phi$
tel que $\overline{U}$ est compact, avec
$M\backslash \Phi$ globalement minimale.
Pour tout $\varepsilon>0$, il existe
une d\'eformation $M^d$ ${\cal C}^{2,\beta}$, $\beta<\alpha$,
avec ${\rm supp} \ d =\overline{U}$, $M^d\supset \Phi$,
$||M^d-M||_{{\cal C}^{2,\beta}}<\varepsilon$,
telle qu'il existe une fonction 
$f^d\in L_{loc}^1(M^d)\cap 
L_{loc,CR}^1(M^d\backslash \Phi)\cap 
{\cal H}({\cal V}(U^d))$ co\"{\i}ncidant
avec $f$ sur $M\backslash U$ et
telle que $||f^d-f||_{L^1}<\varepsilon$.}

\smallskip
La preuve utilise la Proposition 1 et l'in\'egalit\'e
de Carleson sur des d\'eformations successives de $M$ \`a support
la base $U_j\subset M$ de petits wedges ${\cal W}_{A_j,p_j}$
tels que $\bigcup_{j\in J} U_j =M$. Il suffira alors de d\'emontrer
que $L_{loc}^1(M^d) \cap {\cal H}({\cal V}(M^d\backslash \Phi))=
L_{loc,CR}^1(M^d)$. En effet, par
$|\int_M (f^d-f)\overline{\partial}\varphi| \leq C_{\varphi} \varepsilon$
(les mesures sur $M$ et sur $M^d$ \'etant voisines, puisque 
$||M^d-M||_{{\cal C}^{2,\beta}}<\varepsilon$, 
une telle \'ecriture \`a un sens),
pour toute $(m+n,m-1)$-forme \`a support compact,
l'\'egalit\'e $\int_{M^d} f\overline{\partial} \varphi=0$
impliquera $\int_M f\overline{\partial} \varphi=0$, puisque
$\varepsilon$ est arbitraire. $\square$

\smallskip
\noindent
{\it \'Etape 3: ${\cal W}$-\'elimination 
de la singularit\'e.} En utilisant l'hypoth\`ese
<<$M$ et $M\backslash \Phi$ globalement minimales>>, on \'elimine
progressivement les points de $\Phi=N$,
$E$ ou $K$ qui se trouvent \`a l'extr\'emit\'e
d'une courbe int\'egrale par morceaux de $T^cM$ issue
d'un point de $M\backslash \Phi$ et on d\'eforme
ensuite $M$ dans le wedge obtenu au-dessus de chaque
point qui a \'et\'e \'elimin\'e. 
\`A chaque pas, sur une d\'eformation de $M$ encore not\'ee $M$,
la situation se r\'eduit \`a l'\'elimination
d'un seul point $p$ de $\Phi$ dispos\'e comme suit. Il existe
un voisinage $U$ de $p$ dans $M$ et 
$M_1$ une hypersurface ${\cal C}^2$ dans $U$ 
qui partage $U$ en deux composantes
ferm\'ees $M_1^-$ et $M_1^+$, 
$U=M_1^-\cup M_1^+$, $M_1^-\cap M_1^+=M_1$, 
telle que $\Phi \cap U \subset M_1^{-}$,
et il existe un disque $A$ ${\cal C}^{2,\beta}$ 
attach\'e \`a $M_1^+$ avec $A(1)=p$, 
$dA/d\theta (1) \in T_pM_1$. Soient $\omega$ un voisinage de 
$U\backslash \Phi$ dans $\C^{m+n}$ et $f\in {\cal H}(\omega)$. 
Les d\'eformations normales de Tumanov nous permettent de d\'evelopper
$A$ en une famille r\'eguli\`ere en $p$ de disques analytiques 
$A_{s,v}$, $A_{0,0}=A$, qui engendre un wedge ${\cal W}_{A,p}$ :

\smallskip
\noindent
{\sc Lemme 1} (\cite{TU}.) -- {\it Il existe une famille 
r\'eguli\`ere $A_{s,v}$
attach\'ee \`a $(M\cap U)\cup \omega$.} $\square$

\smallskip 
Cependant, le bord de ces disques peut toucher la singularit\'e $\Phi$
(en fait, $A_{s,v}(b\Delta) \cap \Phi \subset A_{s,v}(b\Delta \cap
\Delta_1)$) et le th\'eor\`eme d'approximation de Baouendi-Treves
n'est plus valable.  Heureusement, le principe de continuit\'e et une
propri\'et\'e d'isotopie des disques \`a un point nous permet de
prolonger $f$ \`a ${\cal W}:={\cal W}_{A,p}$ moins un ensemble ${\cal
E}\subset {\cal W}$, {\it i.e.} $F\in {\cal H}(\omega \cup ({\cal W}
\backslash {\cal E}))$, comme suit:

\smallskip
\noindent
{\sc D\'efinition 4.} Un disque plong\'e 
$A: \overline{\Delta}\to \C^{m+n}$
est dit {\it $b$-isotope \`a un point dans $\omega$} s'il existe une 
application ${\cal C}^1$ $[0,1]\times 
\overline{\Delta} \ni (t,\zeta)\mapsto 
A_t(\zeta)\in \C^{m+n}$ telle que $A_t(b\Delta)\subset\omega$,
$A_0=A$, chaque $A_t$ est
un disque analytique plong\'e pour $0\leq t< 1$ et
$A_1$ est une application constante 
$\overline{\Delta} \to \{pt\}\in \omega$.

\smallskip
\noindent
{\sc Lemme 2.} -- {\it Sous les conditions des
th\'eor\`emes 2, 3 et 4, tout disque $A_{s,v}$ tel
que $A_{s,v}(b\Delta\cap \Delta_1) \cap \Phi =\emptyset$
est $b$-isotope \`a un point dans 
$\omega$.} $\square$

\smallskip
\noindent
{\sc Lemme 3.} -- {\it Soit $A: b\Delta \to \omega$, 
$\overline{\Delta} \to \C^{m+n}$ $b$-isotope \`a un point dans 
$\omega$. Alors, pour toute fonction holomorphe $f\in {\cal H}(\omega)$,
il existe $F\in {\cal H}(\omega \cup 
{\cal V}(A(\overline{\Delta})))$
telle que $F\equiv f$ dans ${\cal V}(A(b\Delta))$.} $\square$

\smallskip
L'ensemble
${\cal E}:=\{z_1=A_{s,v}(\zeta_1)\in \C^{m+n}: \ \zeta_1 \in
\stackrel{\circ}{\Delta}_1, A_{s,v}(b\Delta) \cap \Phi\neq \emptyset\}$
o\`u l'on n'a pas prolong\'e
est feuillet\'e par des courbes holomorphes. Gr\^ace aux disques
attach\'es \`a $M$, on se ram\`ene donc \`a \'eliminer
la singularit\'e ${\cal E}\backslash \omega$ pour 
$F\in {\cal H}(\omega \cup ({\cal W} \backslash {\cal E}))=
{\cal H}(\omega \cup({\cal W} \backslash ({\cal E}\backslash
\omega)))$.
\`A moins que ${\rm codim}_{\cal W} {\cal E}=
{\rm codim}_M \Phi$, le bord de presque tout
disque $A_{s,v}$ ne touche en g\'en\'eral 
$\Phi$ que sur un ensemble
de mesure nulle de $b\Delta$.
Le fait que de nombreux disques $A_{s,v}$ satisfont
$A_{s,v}(b\Delta \cap \Delta_1) \not\subset {\cal E}$ pr\`es de
$\zeta=1$, {\it i.e.} que ${\cal E}\backslash \omega\neq {\cal E}$,
est crucial pour la suite.

La structure de ${\cal E}$ d\'epend des cas:

$\bullet$ Si ${\cal H}^{d-3}_{loc}(E)<\infty$, alors
${\cal H}_{loc}^{2m+2n-2}({\cal E})<\infty$. Soit
$f\in L^1(U)$. Dans ce cas, 
$F\in L^1({\cal W})$ et on d\'emontre que ${\cal E}\backslash \omega$ est une
singularit\'e \'eliminable pour $F$ gr\^ace \`a un principe de m\'eromorphie
s\'epar\'ee d\^u \`a Shiffman. En effet, on a:

\smallskip \noindent {\sc Lemme 4.} -- {\it Soit $P\subset\subset
{\cal W}$ un polydisque. Alors pour presque tout disque de coordonn\'es
$D\subset P$, $D\cap {\cal E}$ consiste en un nombre fini de points et
$F|_D$ est m\'eromorphe sur $D$ \`a p\^oles d'ordre au plus un.}

\smallskip
Le lemme s'applique \`a
$F\in L^1({\cal W})\cap {\cal H}(\omega \cup ({\cal W}
\backslash {\cal E}))$, d'o\`u $F$ est m\'eromorphe
dans ${\cal W}$. Si l'ensemble
polaire $P_F$ de $F$ est non vide, il ne peut pas contenir
une courbe $A_{s,v}(\Delta)$ telle que 
$A_{s,v}(\stackrel{\circ}{\Delta}_1)\cap
\omega \neq \emptyset$. Donc $P_F$ est constitu\'e
des disques dont le bord est enti\`erement contenu, pour $\zeta$
pr\`es de $1$, dans $E$. L'ensemble ${\cal E}_1$ correspondant dans
${\cal W}$ satisfait maintenant $H^{2m+2n-3}({\cal E}_1)=0$, mais
alors ${\cal H}({\cal W} \backslash {\cal E}_1) ={\cal H}({\cal W})$,
donc $P_F=\emptyset$. $\square$

\smallskip

$\bullet$ Dans la situation du th\'eor\`eme 4, 
l'\'enonc\'e suivant s'applique \`a (i), (ii) et (iii) :

\smallskip
\noindent
{\sc Proposition 4.} -- {\it Soit $\Lambda$ une hypersurface d'un ouvert 
${\cal U} \subset \C^{m+n}$
et $\Phi$ un ferm\'e de $\Lambda$ qui ne contient pas d'orbite CR de $\Lambda$.
Alors $\forall \ f\in {\cal H}({\cal U}
\backslash \Phi), \exists \ F\in {\cal H} ({\cal U})$ telle que
$F|_{{\cal U} \backslash \Phi} =f$.}

\smallskip
Lorsque ${\rm codim}_M \ N=2$ (Th\'eor\`eme 4), 
${\cal E}$ est une hypersurface de ${\cal W}$ feuillet\'ee
par les $A_{s,v}(\stackrel{\circ}{\Delta}_1)$, mais
tout bord de disque $A_{s,v}(b\Delta\cap \Delta_1)$ ne rencontre
$N$ qu'en au plus un point. Donc ${\cal E}\backslash \omega$ ne 
contient pas d'orbite CR de ${\cal E}$: la proposition 4 
s'applique. $\square$

\smallskip
$\bullet$ Lorsque $M$ est ${\cal C}^{\omega}$ et
$H^{d-2}({\cal E})=0$, ${\cal E}$ satisfait $H^{2m+2n-1}(E)=0$.
Dans ce cas, le feuilletage de ${\cal E}$ est analytique r\'eel et
on applique pr\`es d'un point de
$b\omega\cap {\cal E}$ l'\'enonc\'e (\cite{MP2}):

\smallskip
\noindent
{\sc Proposition 5.} -- {\it Soient ${\cal U}\subset \C^{m+n}$ un ouvert
${\cal C}^{\omega}$-feuillet\'e par 
des courbes holomorphes $A_{\vartheta}, \vartheta \in D$,
${\cal U}=\cup_{\vartheta} A_{\vartheta}$, $0\in {\cal U}$,
$D\subset \R^{2m+2n-2}$ un ouvert, $0\in D$, 
${\cal G}\subset D$ un ferm\'e avec ${\cal H}^{2m+2n-3}({\cal G})=0$ et
$M_1$ une hypersurface ${\cal C}^1$ dans ${\cal U}$, $0\in M_1$,
$T_0M_1+T_0A_0=T_0\C^{m+n}$, et posons ${\cal E}:= (\cup_{\vartheta\in
{\cal G}} A_{\vartheta}) \cap M_1^-$. Alors ${\cal H}({\cal U}
\backslash {\cal E})={\cal H}({\cal U})$.}

\smallskip
Ici, l'argument utilise une complexification 
de courbes r\'eelles du feuilletage
pour avoir la $b$-isotopie. $\square$

\smallskip
\noindent
{\it \'Etape 4: valeurs au bord dans $L^1$.} 
Pour conclure, Il reste \`a estimer l'extension $F$
dans le  wedge ${\cal W}$ en norme $L^1$.
On applique \`a 
$\Phi\subset M_1^-$ l'\'enonc\'e:

\smallskip
\noindent
{\sc Proposition 6.} -- {\it Si 
${\cal H}_{loc}^{d-2}(\Phi)<\infty$ et si
${\cal H}(\omega)$ 
se prolonge holomorphiquement
\`a un wedge en $p\in M_1^-$, alors $p$ est
$L^1$-\'eliminable, {\it i.e.} 
<<$\Phi$ ${\cal W}$-\'eliminable>>
entra\^{\i}ne
<<$\Phi$ $L^1$-\'eliminable>>.}

\smallskip
{\it Preuve:} Les disques $A_{s,v}$ et l'in\'egalit\'e
de Carleson donnent un contr\^ole en norme $L^1$ pour les
disques ne touchant pas $\Phi$
({\it i.e.} presque tout disque, parce que
${\cal H}^{d-2}_{loc}(\Phi) < \infty$),
donc de l'extension dans le wedge.
L'extension appartient enfin \`a ${\rm H}_{\rm a}^1({\cal W})$,
ce qui ach\`eve la preuve. $\square$

\medskip
\noindent
{\it Preuve du Th\'eor\`eme 5.} 
Compte tenu de l'existence des familles r\'eguli\`eres
de disques analytiques et de la Proposition 2,
en suivant la d\'emonstration du Th\'eor\`eme 1 dans 
Kytmanov-Rea \cite{KR}, on obtient le
Th\'eor\`eme 5. $\square$

\end{document}